# Integrated fractional Brownian motion: persistence probabilities and their estimates

**G.Molchan**

Institute of Earthquake Prediction Theory and Mathematical Geophysics,
Russian Academy of Science, Profsoyuznaya 84/32, Moscow, Russia

***Abstract***. The problem is a log-asymptotics of the probability that the Integrated fractional Brownian motion of index 0<H<1 does not exceed a fixed level during long time. For the growing time interval (0,T) the hypothetical log-asymptotics is (H(H-1)+o(1))Log T. In support of the hypothesis, we update our earlier estimates of the probability and give analytical proofs.

## 1. The problem

Let $x(t), x(0)=0$ be a real-valued stochastic process with the polynomial asymptotics:

$$P(x(t)<1, t\in\Delta_T)=T^{-\theta_x+o(1)},\ \ T\to\infty\ ,\quad \Delta_T=T\Delta, \tag{1.1}$$

where $\Delta$ is some bounded interval containing 0. In that case, $\theta_x$ is known as the *persistence* exponent. The problem of exact values of $\theta_x$ is usually a difficult task. An overview and some new results are presented in (Bray et al.,2013; Aurzada and Simon, 2015; Profeta and Simon, 2015; Molchan, 2017;Aurzada et al., 2018).

Below we consider the process which is obtained by integrating the fractional Brownian motion, i.e., $I_H(t)=\int_0^t w_H(s)ds$ where $\mathrm{w}_{\mathrm{H}}(t), t\in R^1$ is a Gaussian process with the correlation function

$$\mathrm{Ew}_{\mathrm{H}}(t)w_H(s)=1/2(|t|^{2H}+|s|^{2H}-|t-s|^{2H})\,.$$

As is known (Molchan 1999),

$$\theta_{w_H}=1-H\ ,\ \Delta_T=(0,T)\,. \tag{1.2}$$

A similar result holds for the persistence exponent of $I_H(t)$, but in a bilateral growing interval:

$$\theta_{I_H}=1-H,\qquad \Delta_T=(-T,T)\,, \tag{1.3}$$

(Molchan , 2017) . In the case $H=1/2$ the trajectories of $I_H(t)$ are independent in the intervals $t<0$ and $t>0$. Therefore, putting $H=1/2$ in (3), we arrive at the Sinai (1992) result for the integrated Brownian motion in $\Delta_T=(0,T)$, namely,

$$\theta_{I_{1/2}}=1/4\ ,\qquad \Delta_T=(0,T)\,. \tag{1.4}$$

In the general case of $H$, the exact value of the persistence exponent $\theta_{I_H}$ for $\Delta_T = (0,T)$ is unknown, and remains an important unsolved problem. Molchan and Khokhlov (2004) analyzed the exponent $\theta_{I_H}$ theoretically and numerically to put forward the following

***Hypothesis***: $\theta_{I_H} = H(1-H)$ for $\Delta_T = (0,T)$.

The hypothetical equality $\theta_{I_H} = \theta_{I_{1-H}}$ is quite interesting, because the exponents are related to processes with fundamentally different probabilistic properties. This difference is well reflected in the result (1.3) for the bilateral growing interval $\Delta_T = (-T,T)$. In support the hypothesis, below we update our earlier estimates of $\theta_{I_H}$ for $\Delta_T = (0,T)$ (see Molchan , 2012) and give analytical proofs for them . The main result is the following

**Proposition 1.** For the integrated fractional Brownian motion in $\Delta_T = (0,T)$

**a)** $\theta_{I_H} \geq \theta_{I_{1-H}}$, $\quad 0 < H \leq 0.5$;

**b)** $0.5(H \wedge \overline{H}) \leq \theta_{I_H} \leq H \wedge \overline{H}$, $\quad \overline{H} = 1-H$;

**c)** $\theta_{I_H} \leq 1/4 \vee \sqrt{(1-H^2)/12}$.

**Remark.** The above estimates are closest to the hypothetical values near the indexes $H = 0,\ 1/2$ , 1 and are well matched with possible symmetry of $\theta_{I_H}$ relative to the point $H = 1/2$.

## 2. Modification of the problem

Any self-similar process $x(t)$ in $\Delta_T$ =(0,T) generates a *dual stationary process* $\tilde{x}(s) = e^{-hs}x(e^s)$, $s < \ln T = \tilde{T}$, where $h$ is the self-similarity index of $x(t)$. For a large class of Gaussian processes, relation (1.1) induces the dual asymptotics

$$P(\tilde{x}(s) \leq 0, 0 < s < \tilde{T}) = \exp(-\tilde{\theta}_x \tilde{T}(1+o(1))), \quad \tilde{T} \to \infty$$

with the same exponent $\tilde{\theta}_x = \theta_x$. This is particularly true for the processes $I_H(t)$ and $w_H(t)$ [Molchan 1999, 2008]. The equality $\tilde{\theta}_x = \theta_x$ reduces the original problem to the estimation of $\tilde{\theta}_x$. Non-negativity of the correlation function of $\tilde{x}(s)$, $\tilde{B}_x(s)$, guarantees the existence of the exponent $\tilde{\theta}_x$ [ Li and Shao, 2004]. This is the case of $\tilde{I}_H(t)$. The inequality of two correlation functions, $\tilde{B}_1(s) \leq \tilde{B}_2(s)$, $\tilde{B}_i(0) = 1$, implies, by Slepian's lemma, [Lifshits,1995], the inverted inequality for the relevant exponents: $\tilde{\theta}_1 \geq \tilde{\theta}_2$. Therefore, Proposition 1 follows from

**Proposition 2.** Let $\tilde{B}_{I_H}(t)$ and $\tilde{B}_{w_H}(t)$ be the correlation functions of the stationary processes that are dual to $I_H$ and $w_H$, respectively. Then

**a)** $\tilde{B}_{I_H}(t) \leq \tilde{B}_{I_{1-H}}(t)$, $\quad 0 < H \leq 0.5$, $\quad$ (2.1)

**b)** $\tilde{B}_{I_H}(t) \ge \tilde{B}_{w_{1-H}}(t)$, $\quad 0 < H \le 0.5$, (2.2)

**c)** $\tilde{B}_{I_{1/2}}(t) \le \tilde{B}_{I_H}(t) \le \tilde{B}_{I_{1/2}}(2(1-H)t)$, $\quad 0.5 < H < 1$, (2.3)

**d)** $\tilde{B}_{I_H}(t) \ge \tilde{B}_{I_{1/2}}(2t\sqrt{(1-H^2)/3})$, $\quad 0.25 \le H \le 0.5$. (2.4)

**The relationship between Propositions 1 and 2.**

Assume that Proposition 2 holds, then

1) (2.1) entails $\tilde{\theta}_{I_H} \ge \tilde{\theta}_{I_{1-H}}$ for $0 < H \le 0.5$ , which supports Proposition 1(**a**);

2) (2.2) and (1.2) entail $\tilde{\theta}_{I_H} \le \tilde{\theta}_{w_{1-H}} = H$ for $0 < H \le 0.5$, which supports Proposition 1(**b**,right) for the case $0 < H \le 0.5$;

3) (2.3, right) and (1.4) entail $\tilde{\theta}_{I_H} \ge p(H)/4$ with $p(H) = 2(1-H)$ for $0.5 < H < 1$ . Here we use the obvious fact that the exponent for $\tilde{x}(pt)$ is $p\tilde{\theta}_x$. Using the inequality $\tilde{\theta}_{I_H} \ge \tilde{\theta}_{I_{1-H}}$ for $0 < H \le 0.5$, we shall also have $\tilde{\theta}_{I_H} \ge p(1-H)/4 = H/2$ for $H \le 0.5$ . Finally, one has $\tilde{\theta}_{I_H} \ge 0.5H \wedge (1-H)$, which supports Proposition 1 (**b**,left);

4) (2.3, left) and (1.4) entail $\tilde{\theta}_{I_H} \le \tilde{\theta}_{I_{1/2}} = 1/4$ for $0.5 < H < 1$, which supports Proposition 1(**c**) for $0.5 < H < 1$.

5) (2.4) and (1.4) entail $\tilde{\theta}_{I_H} \le p(H)/4$ with $p(H) = 2\sqrt{(1-H^2)/3}$, which supports Proposition 1(**c**) for $0.25 \le H \le 0.5$. This estimate is trivial for the interval $0 < H \le 0.25$ because $\tilde{\theta}_{I_H} \le H$ according to Proposition 1(**b**) .

6) It remains to show that $\theta_{I_H} \le 1-H$ (Proposition 1,**b**). This fact obviously follows from (1.3), because $\theta_{I_H}(\Delta_T^1) \le \theta_{I_H}(\Delta_T^2) = 1-H$ for $\Delta_T^1 = (0,T) \subseteq \Delta_T^2 = (-T,T)$. In addition, $\theta_{I_H} = \tilde{\theta}_{I_H}$.

## 3. **Proof of Proposition 2.**

Below we use the dual processes to the integrated fractional Brownian motion $I_H(t)$ and to the fractional Brownian motion $w_H(t)$. They are stationary and have the following correlation functions:

$$\tilde{B}_{I_H}(t) = (2+4H)^{-1}[(4+4H)\cosh(Ht) - 2\cosh((1+H)t) + (2\sinh(t/2))^{2H+2}], \quad (3.1)$$

$$\tilde{B}_{I_{1/2}}(t) = 1/2(3\exp(-|t|/2) - \exp(-3|t|/2)), \quad (3.2)$$

$$\tilde{B}_{w_H}(t) = \cosh(tH) - (2\sinh(t/2))^{2\mathrm{H}}/2. \quad (3.3)$$

**Proof of the relation (2.1):** $\widetilde{B}_{I_H}(t) \le \widetilde{B}_{I_{1-H}}(t)$, $\quad 0 < H \le 0.5$.

***Notation:*** $\quad x = \exp(-t)$, $\alpha = 1-2H$ ; $\bar{x} = 1-x$ .

By (3.1),

$$2(4-\alpha^2)x^{3/2}(\widetilde{B}_{I_H}(t) - \widetilde{B}_{I_{1-H}}(t)) = U(x,\alpha) - U(x,-\alpha) := \Delta(x,\alpha) \tag{3.4}$$

where

$$U(x,\alpha) = (2+\alpha)[\bar{x}^{3-\alpha} - 1 + (3-\alpha)x + (3-\alpha)x^{2-\alpha} - x^{3-\alpha}]x^{\alpha/2}$$

We have to show that $\Delta(x,\alpha) \le 0$ on $S = [0,1]^2$.

***Step 1 :*** $0 \le x \le 0.5,\ 0 \le \alpha \le 1.$

By the Taylor expansions,

$$\bar{x}^{3-\alpha} = 1 - (3-\alpha)x + (3-\alpha)(2-\alpha)x^2/2 - (3-\alpha)(2-\alpha)(1-\alpha)x^3/6 - r(x|\alpha)\alpha x^2/6, \tag{3.5}$$

where

$$r(x|\alpha) = (3-\alpha)(2-\alpha)(1-\alpha)x^2 \int_0^1 (1-u)^3(1-ux)^{-1-\alpha}\,du \ .$$

Since $(1-ux)^q \ge 1$ for $q < 0$,

$$r(x|\alpha) \ge (3-\alpha)(2-\alpha)(1-\alpha)x^2/4 \ , \ |\alpha| < 1 \tag{3.6}$$

and

$$\Delta(x,\alpha) = -\alpha x^2[f(x|\alpha) + R(x|\alpha)]/6, \tag{3.7}$$

where

$$f(x|\alpha) = 3(2-\alpha)(1+\alpha)x^\alpha - (2-\alpha)(5+2\alpha-\alpha^2)x^{1+\alpha}$$
$$+ 3(2+\alpha)(1-\alpha) - (2+\alpha)(5-2\alpha-\alpha^2)x \tag{3.8}$$

and

$$R(x|\alpha) = (2+\alpha)r(x|\alpha) + (2-\alpha)r(x|-\alpha) \ . \tag{3.9}$$

Due to (3.6), we have

$$R(x|\alpha) \ge (4-\alpha^2)(3+\alpha^2)x^2/2 \tag{3.10}$$

Therefore

$$[f(x|\alpha)+R(x|\alpha)]\geq\tilde{f}(x|\alpha)=f(x|\alpha)+(4-\alpha^2)(3+\alpha^2)x^2/2$$

Now we are going to show that $\tilde{f}(x|\alpha)\geq 0$.

• The function $x\to f(x|\alpha)$ is concave, and $f(0|\alpha)>0$. In addition,

$$2f(0.5|\alpha)=(2-\alpha)(1+4\alpha+\alpha^2)2^{-\alpha}+(2+\alpha)(1-4\alpha+\alpha^2)$$

is decreasing because

$$2f'_\alpha(0.5|\alpha)=-(7-4\alpha-3\alpha^2)(1-2^{-\alpha})-(2-\alpha)(1+4\alpha+\alpha^2)2^{-\alpha}\ln 2-8\alpha\leq 0$$

Hence,

$$f(0.5|\alpha)\geq f(0.5|0.65)=0.34\geq 0 \text{ for } \alpha\in[0,0.65].$$

As a result ,

$$\tilde{f}\geq f\geq 0 \text{ for } (x,\alpha)\in[0,0.5]\times[0,0.65].$$

• *Assume that* $\alpha>0.65$. In this case the function $x\to\tilde{f}(x|\alpha)$ is concave. Indeed, the expression

$$\varphi(x,\alpha):=-x^{2-\alpha}(2-\alpha)^{-1}\tilde{f}''_x=\alpha(1+\alpha)[3(1-\alpha)+(5+2\alpha-\alpha^2)x]-(2+\alpha)(3+\alpha^2)x^{2-\alpha}$$

is nonnegative because $\varphi(0,\alpha)\geq 0$ , $x\to\varphi(x,\alpha)$ is concave, and $\varphi(0.5,\alpha)\geq 0$. Only the last property is not obvious and needs to be verified.

Using here and below the symbol $\beta=1-\alpha$, one has

$$2\beta^{-1}\varphi(0.5,\alpha)=4-19\beta+10\beta^2-\beta^3+(3-\beta)(4-2\beta+\beta^2)u(\beta), \qquad (3.11)$$

where

$$u(\beta)=(1-2^{-\beta})/\beta\geq u(0.35)\geq 0.61 \qquad 0\leq\beta\leq 0.35.$$

Hence we can continue (3.11) as follows

$$\geq 11.32-25.1\beta+13.05\beta^2-1.61\beta^3\geq 11.32-25.1\beta\geq 0 \quad , \quad 0\leq\beta\leq 0.35.$$

Consequently, the function $\tilde{f}(x|\alpha)$ is concave.

• Since $\tilde{f}(0|\alpha)>0$ , we will have $\tilde{f}(x|\alpha)\geq 0$ for $0\leq x\leq 0.5$ if $\tilde{f}(0.5|\alpha)>0$. One has

$$8\tilde{f}(0.5|\alpha)=2(2-\alpha)(1+4\alpha+\alpha^2)2^{\beta}+4(2+\alpha)(1-4\alpha-\alpha^2)+(4-\alpha^2)(3-\alpha^2).$$

Due to $2^{\beta}\geq 1$, one has

$$8\tilde{f}(0.5|\alpha) \geq \beta(34 - 11\beta + 2\beta^2 - \beta^3) \geq 0 \ , \qquad 0 \leq \beta \leq 1.$$

The proof of the case $0 \leq x \leq 0.5$ is complete.

***Step 2:*** $0.5 \leq x \leq 1, 0 \leq \alpha \leq 0.6$.

***Notation:*** $\beta = 1 - \alpha$, $y = 1 - x$.

By (3.4),

$$\hat{f}(x|\alpha) := -(1-x)^{-1} x^{\alpha/2} \Delta(x,\alpha)$$

$$= [\alpha^2 + \alpha(2-\alpha)y + (2-\alpha)y^2 - (2+\alpha)y^{2-\alpha}](1-y)^{\alpha} \tag{3.12}$$

$$-\alpha^2 + \alpha(2+\alpha)y - (2+\alpha)y^2 + (2-\alpha)y^{2+\alpha}.$$

Define two functions:

$$\hat{R}(y|\alpha) = [1-(1-y)^{\alpha}]/y = \alpha \sum_{k \geq 0} \frac{(\beta)_k}{(2)_k} y^k \geq 0 \ , \tag{3.13}$$

where $(\beta)_л = \beta(\beta+1)....(\beta+k-1)$, and

$$A(y|\alpha) = (2+\alpha)y^{2-\alpha} - (2-\alpha)y^2 \ . \tag{3.14}$$

We may rewrite (3.12) as follows

$$y^{-1}\hat{f}(x|\alpha) = 4\alpha - [(2+\alpha)y^{1-\alpha} - (2-\alpha)y^{1+\alpha}]_{(1)} - 2\alpha y$$

$$[-\alpha^2 - \alpha(2-\alpha)y + A(y|\alpha)]_{(2)} \hat{R}(y|\alpha) \tag{3.15}$$

where the subscripts at the square brackets are used to number the grouped summands.

- Consider the term $[...]_{(1)} := \psi_1(y|\alpha)$. It is concave, $\psi_1(0|\alpha) = 0$, and $\psi_1(1|\alpha) = 2\alpha$. Hence $y \to \psi_1(y|\alpha)$ increases up to the critical point (the extreme point of the function):

$$y^* = \left[\beta \frac{(2+\alpha)(1-\alpha)}{(2-\alpha)(1+\alpha)}\right]^{1/(2\alpha)} .$$

If $\qquad u(\alpha) := (2-\alpha)(1-\alpha)(1-4^{-\alpha}) - 2\alpha > 0$ ,

then $\qquad y^* \geq 0.5$ and

$$[...]_{(1)} = \psi_1(y|\alpha) \leq \psi(0.5|\alpha) \quad \text{for} \ \ 0 \leq y \leq 0.5. \tag{3.16}$$

The function $u(\alpha)$ is concave, $u(0)=0$, and $u(0.6)=0.065>0$.

Hence (3.16) holds for $0\le\alpha\le 0.6$.

• Consider the second term $[...]_{(2)} := \psi_2(y|\alpha)$ in (3.15):

By (3.14),

$$A(y|\alpha)=(2-\alpha)(y^{2-\alpha}-y^2)+2\alpha y^{2-\alpha}\ge 2\alpha y^{2-\alpha}\ .$$

Therefore,

$$\psi_2(y|\alpha)\ge -\alpha[\alpha+(2-\alpha)y-2y^{2-\alpha}] := -\alpha\varphi(y|\alpha)\,.$$

The function $\varphi(y|\alpha)$ is concave on $0\le y\le 1$, and $\varphi(0|\alpha)=\varphi(1|\alpha)=0$. Therefore $\varphi(y|\alpha)\ge 0$.

The critical point of $y\to\varphi(y|\alpha)$ is $y^\bullet=2^{-1/(1-\alpha)}\le 1/2$. Hence,

$$[...]_{(2)}=\psi_2(y|\alpha)\ge -\alpha\varphi(y^\bullet|\alpha)=-\alpha[\alpha+(1-\alpha)2^{-1/(1-\alpha)}]. \tag{3.17}$$

• We now return to (3.15). Using (3.16), (3.17) , and the simple relations

$$2\alpha y\le\alpha \quad \text{for } y\le 0.5 \quad \text{and} \quad \hat{R}(y|\alpha)\le \hat{R}(0.5|\alpha)\ ,$$

we obtain $\quad (\alpha y)^{-1}\hat{f}(x|\alpha)\ge 3-w(\alpha)\ ,\quad$ where

$$w(\alpha)=[(2+\alpha)2^{\alpha}-(2-\alpha)2^{-\alpha}](2\alpha)^{-1}+2(\alpha+(1-\alpha)2^{-(1-\alpha)^{-1}}(1-2^{-\alpha})$$

$$=2\alpha^{-1}\sinh(\alpha\ln 2)+\cosh(\alpha\ln 2)+2(1-2^{-\alpha})(1-e^{-\lambda})/\lambda\ ,\ \lambda=(1-\alpha)^{-1}\ .$$

We can see that $w(\alpha)$ is an increasing function. Therefore,

$$(\alpha(1-x))^{-1}\hat{f}(x|\alpha)\ge 3-w(\alpha)\ge 3-w(0.6)=0.03>0, \qquad (x,\alpha)\in[0.5,1]\times[0,0.6]\ .$$

In other words, see (3.12), $\Delta(x,\alpha)\le 0$ for $(x,\alpha)\in[0.5,1]\times[0,0.6]$.

***Step 3:*** $0.6\le\alpha\le x\le 1$.

Regroup the elements in (3.15) as follows

$$y^{-1}\hat{f}(x|\alpha)=(2+\alpha)(1-y^{\beta})-(\hat{R}-\alpha)-[(2+\alpha)\beta^2+(3-\beta^2)\alpha y-\alpha(2+\alpha)y^{1+\beta}]_{(1)}$$

$$+[(\beta-y+\beta(1-\beta)+\beta^2 y+(2+\alpha)y^{1+\beta})(\hat{R}(y|\alpha)-\alpha)]_{(2)}\,. \tag{3.18}$$

• By (3.13),

$$\hat{R}-\alpha \geq \alpha\beta y/2 \qquad , \quad \alpha+\beta=1=x+y \tag{3.19}$$

Consider

$$\Phi(y|\beta)=\frac{\hat{R}-\alpha}{\alpha\beta y}=\sum\nolimits_{k\geq 0}\frac{(\beta+1)_k}{(2)_{k+1}}y^k \geq 0\,.$$

This function is increasing in both arguments. Hence

$$\frac{\hat{R}-\alpha}{\alpha\beta y} \leq \Phi(\mathrm{y}_0|\beta_0)=\frac{1-x_0^{\alpha}-\alpha_0 y_0}{\alpha_0\beta_0 y_0}, \quad 0\leq y\leq y_0, \quad 0\leq\beta\leq\beta_0 \quad . \tag{3.20}$$

In particular

$$\hat{R}-\alpha \leq \bar{k}\beta y\,,\ \bar{k}=0.625, 0\leq y,\beta\leq 0.4\,. \tag{3.21}$$

- Consider (3.18). We have $[...]_{(2)}\geq 0$ as a consequence of $\beta\geq y$ and (3.19). Using (3.21) , one has

$$(2+\alpha)^{-1}y^{-1}\hat{f}(x|\alpha)\geq 1-[(1-\alpha y)y^{\beta}+\beta^2+(2+\alpha)^{-1}(\bar{k}\beta+3-\beta^2)\alpha y]:=\varphi(y|\beta)\,.$$

It is easy to show that $\varphi(y\mid\beta)$ is decreasing on $0\leq y\leq 1$. Hence, $\varphi(y\mid\beta)\geq\varphi(\beta\mid\beta)$ for $0\leq y\leq\beta$. The non-negativity of $\varphi(\beta\mid\beta)$ is equivalent to the relation

$$u(\beta):=[(3-4\beta+3\beta^2+(1-\beta)\beta^2]v(\beta)-(\bar{k}+4)+(\bar{k}+3)\beta-\beta^2\geq 0\,. \tag{3.22}$$

Here $v(\beta)=(1-\beta^{\beta})/\beta^2$ is convex on $0\leq\beta\leq 1$ because

$$\beta^4 v''(\beta)=6-\beta^{\beta}[(\beta\ln\beta+\beta-2)^2-\beta-2]\geq 8+\beta-(\beta\ln\beta+\beta-2)^2\geq 8+\beta-(2+e^{-2})\geq 0\,.$$

Therefore,

$$v(\beta)\geq v(0.4)+v'(0.4)(\beta-0.4)=5.898-9.952\beta\ .$$

Using this inequality, we can continue (3.22) as follows

$$u(\beta)\geq 13.069-49.823\beta+56.502\beta^2-29.856\beta^3$$

$$\geq 13-50\beta+26\beta^2\geq 0\,.$$

Here we use the $\bar{k}$ -value from (3.21).The final inequality is obvious for $0\leq\beta\leq 0.4$.

Hence $\hat{f}(x|\alpha)\geq 0$. The proof of the case $0.6\leq\alpha\leq x\leq 1$ is complete.

***Step 4:*** $0.5\leq x\leq\alpha, \alpha\geq 0.6$ or $0.5\geq y\geq\beta, \beta\leq 0.4$.

Regroup the elements in (3.18) again

$$y^{-1}\hat{f}(x|\alpha) = (2+\alpha)(1-y^{\beta})x^{\alpha} + \beta(1+\alpha x)(\hat{R}-\alpha) +$$

$$(1+\alpha)\beta^2 y - (1-2y)(\hat{R}-\alpha) - \beta[(2+\alpha)\beta + y] \ . \quad (3.23)$$

• Due to (3.19), (3.20), we have

$$\hat{R}-\alpha = K\alpha\beta y,\ K \in [\underline{k}, \overline{k}],$$

where

$$\overline{k} = \Phi(0.5\,|\,0.4) = 0.681 \text{ for } 0 \le \beta \le 0.4\,, 0 \le y \le 0.5\,,$$

$$\underline{k} = \Phi(y_0\,|\,\beta_0) \qquad \text{for} \qquad \beta \ge \beta_0, y \ge y_0\,. \quad (3.24)$$

Therefore

$$y^{-1}\hat{f}(x|\alpha) \ge [(2+\alpha)(1-y^{\beta})x^{\alpha}]_{(1)} - [((1-2y)\overline{k}+1)y\beta + (2+\alpha)\beta^2]_{(2)} + m\beta^2 y\,, \quad (3.25)$$

where

$$m = (1+\alpha x)\alpha\underline{k} + \alpha + 1 \ .$$

Note that

$$m \ge \underline{m} = 5\underline{k}/8 + 1.5 \quad \text{for } \alpha \ge 0.5, x \ge 0.5\,.$$

The term $[...]_{(1)}$ in (3.25) is a decreasing function of $y$. The same is true for $(-1)[...]_{(2)}$ in the $y-$ interval $I = [0, (\overline{k}+1)/4\overline{k}]$. We have $\overline{k} \le 1$ in the case under consideration. Therefore, $I \supset [0,0.5]$.

To estimate the third term in (3.25), recall that $\beta \le y$ and $\underline{k} \ge \Phi(0|0) = 0.5$. Hence

$$m\beta^2 y \ge \underline{m}\beta^3, \qquad \underline{m} = 5/16 + 1.5 = 1.8125\,.$$

As a result,

$$y^{-1}\hat{f}(x|\alpha) \ge \varphi(y_0|\beta), \qquad \beta \le 0.4, \beta \le y \le y_0 \le 0.5\,, \quad (3.26)$$

where

$$\varphi(y|\beta) = [(2+\alpha)(1-y^{\beta})x^{\alpha}]_{(1)} - [((1-2y)\overline{k}+1)y\beta + (2+\alpha)\beta^2]_{(2)} + m\beta^2 y \quad (3.27)$$

• For $y_0 = 0.5$, we have $\varphi(y_0|\beta) \ge 0$ if

$$\psi(\beta) := (3-\beta)[(2^{\beta}-1)/\beta - 2\beta] - 1 + 2\underline{m}\beta^2 \ge 0\,. \quad (3.28)$$

But

$$(2^{\beta}-1)/\beta \geq \ln 2 + \beta(\ln 2)^2/2 .$$

Therefore,

$$\psi(\beta) \geq 1.079 - 5.97\beta + 5.38\beta^2 \geq 0 \quad , \; 0 \leq \beta \leq 0.2 .$$

Hence, $\hat{f} \geq 0$ for $\beta \leq 0.2, \beta \leq y \leq 0.5$.

• For $y_0 = 0.4$, $\varphi(y_0 \mid \beta) \geq 0$ if

$$\beta^{-1}(1-y_0{}^{\beta})x_0{}^{\alpha} \geq \beta + (3-\beta)^{-1}(a - \underline{m}\beta^2), \qquad x_0 = 1 - y_0, \tag{3.29}$$

where $a = (1 - 2y_0\bar{k} + 1)y_0 = 0.4545$, $\bar{k} = \Phi(0.5 \mid 0.4) = 0.681$, $\underline{m} = 1.8125$.

To prove (3.29), note that the left part of (3.29) can be represented as

$$\beta^{-1}(1-y_0{}^{\beta})x_0{}^{\alpha} = x_0 \int_{\ln(x_0/y_0}^{\ln(1/x_0} e^{t\beta} dt := u(\beta) .$$

Since $u(\beta)$ is convex, $u(\beta) \geq u(0) + u'(0)\beta$ and (3.29) will be true if

$$u(0) + u'(0)\beta \geq \beta + (3-\beta)^{-1}(a - \underline{m}\beta^2),$$

or $$1.649 - 1.463\beta + 1.783\beta^2 \geq 0 .$$

The last relation is obviously true for $\beta \leq 0.4$ .

Consequently, we have proved the case: $\beta \leq y \leq 0.4$ , $\beta \leq 0.4$.

• Consider the last case: $0.4 \leq y \leq 0.5, 0.2 \leq \beta \leq 0.4$ .Then $\underline{k} = \Phi(0.4 \mid 0.2)$ =0.6039 and $m = (1+\alpha x)\alpha\underline{k} + \alpha + 1 \geq 1.9925 =: \underline{m}$ . We need to verify (3.28) for $0.2 \leq \beta \leq 0.4$ using the new parameter $\underline{m}$ . This can be done as above.

## Proof of Relation (2.2):

$$\widetilde{B}_{I_H}(t) \geq \widetilde{B}_{w_{1-H}}(t) \; , \quad 0 < H \leq 0.5 .$$

To prove this relation, consider

$$\Delta(t,H) = (2+2H)\exp(-(1+H)(\widetilde{B}_{I_H}(t) - \widetilde{B}_{w_{1-H}}(t)) .$$

Using (3.1, 3.3) and the notation: $x = \exp(-t)$, $\alpha = 2H$; $\bar{x} = 1 - x$ , we can rewrite our problem as follows:

$$\widetilde{\Delta}(x,\alpha) = \bar{x}^{2+\alpha} - 1 + (2+\alpha)x + (2+\alpha)x^{1+\alpha} - x^{2+\alpha}$$

$$+ (1+\alpha)\bar{x}^{2-\alpha}x^{\alpha} - (1+\alpha)x^{\alpha} - (1+\alpha)x^2 \geq 0, \qquad 0 \leq \alpha, x \leq 1 .$$

One has

$$\widetilde{\Delta}(x,\alpha)/\bar{x} = \bar{x}^{1+\alpha} - 1 + (1+\alpha)x - x^{\alpha}(\bar{x}^{\alpha} + \alpha) + (1+\alpha)\bar{x}^{1-\alpha}x^{\alpha}$$

$$= (1-\bar{x}^{\alpha})^2\bar{x}^{1-\alpha} + \alpha x^{\alpha}(\bar{x}^{1-\alpha} - \bar{x}) + (1-x^{\alpha})(1-(1-\alpha)x - \bar{x}^{1-\alpha}).$$

Obviously

$$(1-(1-\alpha)x - \bar{x}^{1-\alpha}) = x^2\int_0^1 (1-v)(1-vx)^{-\alpha-1}dv > 0.$$

Therefore $\Delta(t,H)$ is nonnegative.

## Proof of Relation (2.3, left):

$$\widetilde{B}_{I_{1/2}} \le \widetilde{B}_{I_H}(t), \quad 0.5 < H < 1,$$

***Notation***: $x = \exp(-t)$; $\alpha = 2H-1, 0 \le \alpha \le 1$.

Consider

$$\Delta(t,H) = (2+4H)[\widetilde{B}_{I_H}(t) - \widetilde{B}_{I_{1/2}}(t)]\exp(-(1+H)t).$$

Due to (3.1, 3.2), the problem looks in terms of the new variables $(x,\alpha) \in S = (0,1)\times(0,1)$ as follows:

$$\widetilde{\Delta}(x,\alpha) = (1-x)^{\alpha+3} - 1 + (3+\alpha)x + (3+\alpha)x^{\alpha+2} - x^{\alpha+3} - (\alpha+2)(3-x)x^{2+\alpha/2} \ge 0.$$

Since $\widetilde{\Delta}(x,0) = 0$, we shall have

$$\widetilde{\Delta}(x,\alpha) = \widetilde{\Delta}(x,\alpha) - \widetilde{\Delta}(x,0)$$

$$= (1-x)^2[(1-x)^{\alpha+1} - (1-x) + \alpha x^{1+\alpha/2}] + \alpha x(1-x^{\alpha/2})(1-x^{1+\alpha/2}) + (3-x)x^2(1-x^{\alpha/2})^2.$$

Using the inequality

$$(1-x)^{\alpha+1} - (1-x) \ge -\alpha x + (1+\alpha)\alpha x^2/2,$$

we can continue

$$\widetilde{\Delta}(x,\alpha) \ge (1-x)^2\alpha x(-1 + (1+\alpha)x/2 + x^{\alpha/2}) + \alpha x(1-x^{\alpha/2})(1-x^{1+\alpha/2})$$

$$= \alpha x^2[(1+\alpha)(1-x)^2/2 + (1-x^{\alpha/2})(2-x-x^{\alpha/2})] \ge 0.$$

## Proof of Relation (2.3,right):

$$\widetilde{B}_{I_H}(t) \le \widetilde{B}_{I_{1/2}}(pt), \; H \ge 1/2, \; p = 2(1-H).$$

***Notation***: $x = \exp(-t)$; $\alpha = 2H-1, 0 \le \alpha \le 1; \beta = 1-\alpha, y = 1-x$.

Consider

$$\Delta(t,H)=(2+4H)[\widetilde{B}_{I_H}(t)-\widetilde{B}_{I_{1/2}}(2(1-H)t)]\exp(-(1+H)t).$$

Due to (3.1, 3.2), the problem looks in terms of the new variables $(x,\alpha)\in S=(0,1)\times(0,1)$ as follows:

$$\widetilde{\Delta}(x,\alpha)=(3+\alpha)(x+x^{\alpha+2})-1-x^{\alpha+3}+(1-x)^{\alpha+3}-3(\alpha+2)x^2+(\alpha+2)x^{3-\alpha}\geq 0\ .$$

By simple algebra one gets

$$-(xy)^{-2}\widetilde{\Delta}(x,\alpha)=[(2+\alpha)\hat{R}(y\mid\beta)+1]\hat{R}(y\mid\alpha)-I(x\mid\alpha),$$

where

$$\hat{R}(y|\alpha)=[1-(1-y)^{\alpha}]/y=\alpha\sum_{k\geq 0}\frac{(\beta)_k}{(2)_k}y^k\geq 0$$

$$I(y\mid\alpha)=[y^{1+\alpha}-1+(1+\alpha)x]x^{-2}=\alpha(1+\alpha)\int_0^1(1-v)(1-xv)^{\alpha-1}dv \tag{3.30}$$

Since $\hat{R}(y\mid a)\geq a$ (see 3.13)), we get

$$-\alpha^{-1}(xy)^{-2}\widetilde{\Delta}(x,\alpha)\geq[(2+\alpha)\beta+1]-I(y\mid\alpha)/\alpha\ . \tag{3.31}$$

The integral representation of $I(y\mid\alpha)$ leads to the following relation

$$I(y\mid\alpha)\leq 0.5\alpha(1+\alpha)(1-x_0)^{\alpha-1},\ 0\leq x\leq x_0$$

and, as a result, to the conclusion:

the right part of (3.31) is non- negative for $0\leq x\leq x_0$ =5/6, if

$$6^{-\alpha}\leq(1+\alpha)^{-1}-\alpha/3\ . \tag{3.32}$$

Using the piecewise linear approximation of the convex function $6^{-\alpha}$ with interpolation between the points 0, 0.2,0.5, and 1, we can easily verify (3.32). This proves the case $0\leq x\leq 5/6, 0\leq\alpha\leq 1$.

Consider the case $5/6\leq x\leq 1$. By (3.30) ,

$$I(y\mid\alpha)=[y^{1+\alpha}-(1+\alpha)y+\alpha]x^{-2}\leq\alpha(6/5)^2.$$

Hence

$$-\alpha^{-1}(xy)^{-2}\widetilde{\Delta}(x,\alpha)\geq[(2+\alpha)\beta+1]-36/25\geq 0, 0\leq\alpha\leq 0.845.$$

Suppose now that $\alpha\geq 0.84$. Combining (3.31) with (3.30), one has

$$-\beta^{-1}y^{-2}\widetilde{\Delta}(x,\alpha)\geq x^2(\alpha(2+\alpha)-1)+1-[y^2(y^{-\beta}-1)\beta^{-1}+y]_{(1)}\ . \tag{3.33}$$

The function $u(y,\beta) := [...]_{(1)}$ is increasing as a function of two variables. This property for the variable x follows from the inequality:

$$u'_y(y,\beta) = 2(y^{1-\beta} - y)\beta^{-1} + (1 - y^{1-\beta}) \ge 0\,.$$

Hence

$$u(y,\beta) \le u(1/6, 0.16) = 0.90975<1,$$

i.e., (3.33) is positive for $\alpha \ge 0.84, x \ge 5/6$ .The proof is complete.

## Proof of relation (2.4):

$$\widetilde{B}_{I_H}(t) \ge \widetilde{B}_{I_{1/2}}(pt)\ ,\quad 1/4 \le H \le 1/2,\quad p = 2\sqrt{(1-H^2)/3}\quad .$$

Let $$\Delta(t,H) = (2+4H)(\widetilde{B}_{I_H}(t) - \widetilde{B}_{I_{1/2}}(pt))e^{-(1+H)t}\,.$$

As a function of $x = \exp(-t)$ and $\alpha = 2H$ , it is

$$\widetilde{\Delta}(x,\alpha) = (2+\alpha)(x + x^{\alpha+1}) - 1 - x^{\alpha+2} + (1-x)^{\alpha+2}$$

$$-3(\alpha+1)x^{1+(\alpha+p)/2} + (\alpha+1)x^{1+(\alpha+3p)/2}\quad . \tag{3.34}$$

We need to verify that $\widetilde{\Delta}(x,\alpha) \ge 0,\ (x,\alpha) \in S = (0,1)\times(1/2,1)$ .

It is easy to check that

$$\widetilde{\Delta}(1,\alpha) = \partial/\partial x\ \widetilde{\Delta}(1,\alpha) = 0,\ \text{and}\ \ (\partial/\partial x)^2\ \widetilde{\Delta}(1,\alpha) = (1+\alpha)(\alpha^2 + 3p^2 - 4)\,.$$

Hence, given $\alpha^2 + 3p^2 = 4$ ,

$$\widetilde{\Delta}(x,\alpha) = \widetilde{\Delta}(x,\alpha) - \widetilde{\Delta}(1,\alpha) - \partial/\partial x\widetilde{\Delta}(1,\alpha)(x-1) - (\partial/\partial x)^2\widetilde{\Delta}(1,\alpha)(x-1)^2/2\,. \tag{3.35}$$

Using notation $y = 1 - x$ and the relation

$$(1-y)^\gamma = 1 - \gamma y + \gamma(\gamma-1)y^2/2 + \gamma(\gamma-1)(\gamma-2)y^3 r(x,\alpha)\,, \tag{3.36}$$

where

$$r(y.\alpha) = 1/2\cdot\int_0^1 (1-u)^2(1-uy)^{\gamma-3}du$$

and

$$1/6 \le r(y,\alpha) \le 1/(2\gamma)\,,\qquad 0 \le \gamma \le 3 \tag{3.37}$$

By (3.34) and (3.35),

$$\widetilde{\Delta}(x,\alpha) = (2+\alpha)(1+\alpha)\alpha(1-\alpha)y^3 r(y,1+\alpha) + (2+\alpha)(1+\alpha)\alpha y^3 r(y,2+\alpha)$$

$$-3(1+\alpha)(1+(\alpha+p)/2)[(\alpha+p)/2](1-(\alpha+p)/2)y^3 r(y,1+(\alpha+p)/2)$$

$$+(1-y)(1+\alpha)[(\alpha+3p)/2]((\alpha+3p)/2)-1)(2-(\alpha+3p)/2)y^3 r(y,(\alpha+3p)/2)$$

$$+y^{2+\alpha} \tag{3.38}$$

The curve $\alpha \to p : \alpha^2 + 3p^2 = 4$ contains the point $(\alpha, p) = (1,1)$ and lies below the tangent $\alpha + 3p = 4$ at that point. Moreover,

$$1 < 3p(0)/2 \le (\alpha + 3p(\alpha))/2 \le 2 \ , \ \ 0 \le \alpha \le 1 \ .$$

Hence, the last two summands in (3.38) are non-negative. Neglecting these terms, we will have

$$(1-x)^{-3}(1+\alpha)^{-1}\widetilde{\Delta}(x,\alpha) \ge (2+\alpha)\alpha(1-\alpha)/6 + (2+\alpha)\alpha/6 - 3(\alpha+p)(2-\alpha-p)/8 := V(\alpha).$$

Since

$$3(\alpha+p)^2 = 2\alpha^2 + 6\alpha p + 4 \ \text{, and } 3p \le 4-\alpha \ ,$$

one has

$$6V(\alpha) = (4-\alpha^2)\alpha - (9/2)(\alpha+p) + (3/2)(\alpha^2 + 3\alpha p + 2)$$

$$\ge (4-\alpha^2)\alpha - (3/2)(1-\alpha)(4-\alpha) - (9/2)\alpha + (3/2)(\alpha^2+2) = (1-\alpha^2)\alpha + 3(2\alpha-1)$$

Hence $\widetilde{\Delta}(x,\alpha) \ge 0$ for $(x,\alpha) \in (0,1)\times(1/2,1)$.

## References.